\documentclass[12pt]{scarticle}

\usepackage{
amsmath,
amsfonts,
amssymb,
amscd,
amsthm, 
calc, 
enumerate, 
}

\usepackage[all]{xy}

\newtheorem{theorem}{Theorem}[section]
\newtheorem{problem}[theorem]{Problem}
\newtheorem{lemma}[theorem]{Lemma}
\newtheorem{corollary}[theorem]{Corollary}
\newtheorem{proposition}[theorem]{Proposition}

\theoremstyle{definition}

\newenvironment{acknowledgments}{\noindent\bf Acknowledgments.\rm}{\nolinebreak\par}

\newtheorem{texample}[theorem]{Example}
\newenvironment{example}{\begin{texample}\rm}
{\nolinebreak\mbox{}\nolinebreak\hfill\raisebox{-.262mm}{$\square$}\end{texample}}

\renewcommand{\leq}{\leqslant}
\renewcommand{\geq}{\geqslant}

\newcommand{\restr}{\mbox{\raisebox{.5mm}{$\upharpoonright$}}}
\newcommand{\smallrestr}{\mbox{\scriptsize $\restr$}}

\newcommand{\concat}{\hspace*{.2mm}\widehat{\mbox{ }}\hspace*{.2mm}}

\newcommand{\cmp}[1]{\overline #1 }
\newcommand{\bigset}[1]{\big\{ #1 \big\}}
\newcommand{\T}{\makebox[4mm][c]{$|_T$}}

\newcommand{\M}{{\mathfrak M}}  
\newcommand{\A}{\mathcal{A}}  
\newcommand{\B}{\mathcal{B}}  
\renewcommand{\S}{\mathcal{S}}   
\newcommand{\C}{\mathcal{C}} 
\newcommand{\D}{\mathcal{D}} 
 
\newcommand{\F}{\mathcal{F}}  
  
\newcommand{\TT}{\mathcal{T}}  
\renewcommand{\P}{\mathcal{P}}  
\newcommand{\X}{\mathcal{X}}  
\newcommand{\bfA}{{\bf A}}
\newcommand{\bfB}{{\bf B}}

\newcommand{\bfS}{{\bf S}}

\newcommand{\0}{{\bf 0}}
\newcommand{\Th}{{\rm Th}}
\newcommand{\join}{+}
\newcommand{\meet}{\times}

\newcommand{\Meet}{{\textstyle \prod}}

\newcommand{\dual}[1]{{\rm dual}(#1)}

\newenvironment{reqlist}{
\begin{list}{-}
 {
\setlength{\parskip}{0cm}
\setlength{\topsep}{.3cm}
\setlength{\partopsep}{0mm}
\setlength{\rightmargin}{0cm}
\setlength{\listparindent}{0cm}
\setlength{\itemindent}{0cm}
\setlength{\parsep}{0mm}
\setlength{\leftmargin}{2.0cm}
\setlength{\labelsep}{1.0cm}
\setlength{\itemsep}{.35cm}
\settowidth{\labelwidth}{zzzzzzzz}
\setlength{\leftmargin}{\labelwidth+\labelsep}
 }
}
{\end{list}}

\begin{document}

\title{On the structure of the Medvedev lattice}

\author{Sebastiaan A. Terwijn\thanks{Institute
of Discrete Mathematics and Geometry,
Technical University of Vienna,
Wiedner Hauptstrasse 8--10/E104,
A-1040 Vienna, Austria,
terwijn@logic.at.
Supported by the Austrian Science Fund FWF
under grant P18713-N12.}
}

\date{\normalsize \today}

\maketitle

\begin{abstract} \noindent
We investigate the structure of the Medvedev lattice as a partial order. 
We prove that every interval in the lattice is either finite, in which 
case it is isomorphic to a finite Boolean algebra, or contains 
an antichain of size $2^{2^{\aleph_0}}$, the size of the lattice itself. 
We also prove that it is consistent that the lattice has chains 
of size $2^{2^{\aleph_0}}$, and in fact that these big chains occur in 
every interval that has a big antichain. 
We also study embeddings of lattices and algebras. 
We show that large Boolean algebras can be embedded into the 
Medvedev lattice as upper semilattices, but that 
a Boolean algebra can be embedded as a lattice only if it is countable. 
Finally we discuss which of these results hold for the closely related 
Muchnik lattice. 
\end{abstract}

%

\section{Introduction} \label{intro} 

Medvedev~\cite{Medvedev} originally introduced the lattice that now bears 
his name in order to establish a connection with intuitionistic logic, 
following up on a rather informal idea of Kolmogorov. 
Later the lattice, which we will denote by $\M$, was studied also 
as a structure of independent interest, being a generalization of 
structures such as the Turing degrees and the enumeration degrees
that are contained in $\M$ as substructures.
For example, Muchnik phrased his original solution to Posts 
problem~\cite{Muchnik1956} as a result in the context 
of the Medvedev lattice. 

Let us briefly recall the definition of $\M$. 
Let $\omega$ denote the natural numbers and let $\omega^\omega$ be the 
set of all functions from $\omega$ to $\omega$ (Baire space). 
A {\em mass problem\/} is a subset of $\omega^\omega$. 
Every mass problem is associated with the ``problem'' of 
producing an element of it. 
A mass problem $\A$ {\em Medvedev reduces\/} to 
mass problem $\B$, denoted $\A\leq_M \B$,  
if there is a partial computable functional  
$\Psi:\omega^\omega\rightarrow\omega^\omega$ defined on all of $\B$
such that $\Psi(\B)\subseteq \A$.
That is, $\Psi$ is a uniformly effective method for     
transforming solutions to $\B$ into solutions to~$\A$.   
The relation $\leq_M$ induces an equivalence relation on   
mass problems: $\A\equiv_M\B$ if $\A\leq_M\B$ and $\B\leq_M\A$. 
The equivalence class of $\A$ is denoted by $\deg_M(\A)$ and is 
called the {\em Medvedev degree\/} of~$\A$. 
We denote Medvedev degrees by boldface symbols. 
There is a smallest Medvedev degree, denoted by~$\0$, 
namely the degree of any mass problem containing a 
computable function, and there is a largest degree~$\bf 1$, 
the degree of the empty mass problem, of which it is 
absolutely impossible to produce an element. 
Finally, it is possible to define a meet operator $\meet$ 
and a join operator $\join$ on mass problems:
For functions $f$ and $g$, as usual define the function $f\oplus g$
by $f\oplus g(2x)=f(x)$ and $f\oplus g(2x+1)= g(x)$. 
Let $n\concat \A = \{ n\concat f: f\in \A \}$, where $\concat$ 
denotes concatenation. 
Define 
$$
\A \join\B = \big \{ f\oplus g: f\in\A \wedge g\in B\big \}
$$
and 
$$
\A \meet\B = 0 \,\widehat{\mbox{ }}\A \cup 1\concat\B. 
$$
The structure $\M$ of all Medvedev degrees, ordered by $\leq_M$ 
and together with $\join$ and $\meet$ is a distributive lattice.  
Medvedev~\cite{Medvedev} also showed that it is possible to define 
an implication operator $\rightarrow$ on $\M$, that is, 
$\M$ is a {\em Brouwer algebra\/}. 
But this will not concern us in the present paper since we will 
mainly be studying $\M$ as a {\em partial order}, 
although the lattice operators on $\M$ will play an important 
role throughout. 
For more information and discussion we refer to the 
following literature.
An early reference is Rogers' textbook \cite{Rogers}, which 
contains a discussion of the elementary properties of $\M$.
Sorbi~\cite{Sorbi} is a general survey paper about $\M$.
Sorbi and Terwijn~\cite{SorbiTerwijn} is a recent paper
discussing the connections with constructive logic.
It also contains an alternative proof of Skvortsova's result
that intuitionistic propositional logic can be obtained as
the theory of a factor of $\M$.
Simpson~\cite{Simpson} surveys
Medvedev reducibility on $\Pi^0_1$ classes, especially with
an eye to the connection with algorithmic randomness.
Binns and Simpson~\cite{BinnsSimpson} are concerned with lattice 
embeddings into the Medvedev and Muchnik lattices of $\Pi^0_1$ classes.

Our notation is mostly standard and follows Odifreddi~\cite{Odifreddi} 
and Kunen~\cite{Kunen}. 
$\Phi_e$ is the $e$-th partial computable functional. 
For $f\in\omega^\omega$ we let $f^-$ be the function with  
$f^-(x) = f(x+1)$ (i.e.\ $f$ with its first element chopped off)    
and for a set $\X\subseteq \omega^\omega$ we let 
$\X^- = \{f^-: f\in\X\}$.     
We use $\mathfrak{2}^n$ to denote the Boolean algebra of all subsets 
of $\{0,\ldots,n-1\}$ under inclusion. 
For countable sets $I\subseteq \omega$ and mass problems 
$\A_i$, $i\in I$, we have the meet operator 
$$
\Meet_{i\in I} \A_i = \bigset{i\,\concat f: i\in I \wedge f\in\A_i}. 
$$
One easily checks that for finite $I$ this is M-equivalent to an iteration 
of the meet operator $\meet$. 
If $a\leq b$ in some partial order, we use the interval notation 
$[a,b] = \{x: a\leq x\leq b\}$. Similarly $(a,b)$ denotes an interval 
without endpoints.

\section{Intervals in the Medvedev lattice} \label{intervals} 

In this section we prove that every interval in the Medvedev lattice is 
either finite of exponential size or contains an antichain of the cardinality 
of the full lattice, namely $2^{2^{\aleph_0}}$. 
We first repeat from Sorbi and Terwijn \cite{SorbiTerwijn} 
the basic strategy for obtaining incomparable elements that avoid upper cones. 
We subsequently generalize this construction to obtain larger 
and larger antichains.

\begin{lemma} {\rm (Sorbi and Terwijn \cite{SorbiTerwijn})}\label{construction}
Let $\A$ and $\B$ be mass problems such that 
\begin{equation} \label{condition} 
\forall \C \subseteq \A \mbox{ finite } \; ( \B\meet \C \not\leq_M \A ). 
\end{equation}
Then there exists a pair $\C_0$, $\C_1$ of M-incomparable mass 
problems $\C_0$, $\C_1 \geq_M \A$ such that 
$\B \meet \C_0$ and $\B\meet\C_1$ are M-incomparable. 
(In particular, neither of $\C_0$ and $\C_1$ is above $\B$.)  
\end{lemma}
\begin{proof}
We want to build $\C_0$ and $\C_1$ above $\A$ in a construction that
meets the following requirements for all $e\in\omega$:

\begin{reqlist}

\item[$R^0_e:$] $\Phi_e(\C_0)\not\subseteq \B\meet\C_1$.

\item[$R^1_e:$] $\Phi_e(\C_1)\not\subseteq \B\meet\C_0$.

\end{reqlist}
The $\C_i\subseteq\A\meet \A\equiv_M \A$ will be built as unions of finite sets
$\bigcup_s \C_{i,s}$, such that $\C_{i,s}\subseteq \A\times\A$ for each pair
$i,s$. We start the construction with $\C_{i,0}=\emptyset$. The idea to meet 
$R^0_e$ is simple: By condition (\ref{condition}) we have at stage $s$ of the 
construction that $\B\meet \C_{1,s} \not\leq_M \A$, so there is a witness 
$f\in\A$ such that $\Phi_e(f)\notin \B\meet \C_{1,s}$. (Either by being 
undefined or by not being an element of $\B\meet \C_{1,s}$.) 
We put such a witness into $\C_0$.
Now this $f$ will be a witness to $\Phi_e(\C_0) \not\subseteq \B \meet \C_1$
provided that we can keep future elements of $1\concat\C_1$ distinct from 
$\Phi_e(f)$. The problem is that some requirement $R^1_i$ may want to put 
$\Phi_e(f)$ into $1\concat \C_1$ because $\Phi_e(f)(0)=1$ and the function
$\Phi_e(f)^- = \lambda x.\,\Phi_e(f)(x+1)$ is the {\em only\/} witness that 
$\Phi_i(\A) \not\subseteq \B \meet \C_0$. To resolve this conflict it suffices 
to complicate the construction somewhat by prefixing all elements of $\A$ by an 
extra bit $x\in\{0,1\}$, that is, to work with $\A\meet\A$ rather than $\A$.
This basically gives us {\em two\/} versions of every potential witness, and we
can argue that either choice of them will be sufficient to meet our needs, so 
that we can always keep them apart. We now give the construction in technical 
detail.

Recall that $f^-(x) = f(x+1)$ and that 
$\X^- = \{f^-: f\in\X\}$.
We build $\C_0$, $\C_1\subseteq \A\meet \A$.

{\em Stage s=0.} Let $\C_{0,0} = \C_{1,0} = \emptyset$.

{\em Stage s+1=2e+1.}
We take care of $R^0_e$.
We claim that there is an $f\in\A-\C^-_{0,s}$ and an $x\in\{0,1\}$ such that
\begin{eqnarray} \label{allesistgut}
\exists h \in \C_{0,s}\cup\{x\concat f\} \big( \Phi_e(h)\notin \B\meet 
(\C^-_{1,s} \meet \C^-_{1,s}) \big).
\end{eqnarray}
Namely, otherwise we would have that for all $f\in\A-\C^-_{0,s}$ and $x\in\{0,1\}$
\begin{eqnarray} \label{supposenot}
\forall h\in\C_{0,s}\cup\{x\concat f\} \big( \Phi_e(h)\in \B\meet (\C^-_{1,s} 
\meet \C^-_{1,s}) \big).
\end{eqnarray}
But then it follows that $\A \geq_M \B\meet (\C^-_{1,s} \meet \C^-_{1,s})$,
contradicting the assumption (\ref{condition}).
To see this, assume (\ref{supposenot}) and let
$$
\D = \C_{0,s}\cup \bigset{x\concat f: x\in\{0,1\} \wedge f\in\A- \C^-_{0,s}}.
$$
Then $\B\meet (\C^-_{1,s} \meet \C^-_{1,s}) \leq_M \D$ via $\Phi_e$. But we 
also have $\D\leq_M \A$, so we have $\B\meet \C^-_{1,s} \leq_M \A$, 
contradicting (\ref{condition}). To show 
that  $\D\leq_M \A$, let $\C_{0,s}^- = \{f_1,\ldots,f_{s}\}$
and let $\tilde{f}_i$, $1\leq i \leq s$,  be finite initial segments such that
the only element of $\C_{0,s}^-$ extending $\tilde{f}_i$ is $f_i$.
(Note that such finite initial segments exist since $\C_{0,s}^-$ is finite.)
Let $x_i$ be such that $x_i\concat f_i \in \C_{0,s}$.
Then $\D\leq_M \A$ via
$$
\Phi(f) =
\begin{cases}
x_i\concat f &  \text{if } \exists i \; \tilde{f}_i \sqsubseteq f, \\
0\concat f   &  \text{otherwise.}
\end{cases}
$$

So we can choose $h$ as in (\ref{allesistgut}).
Put $h$ into $\C_{0,s+1}$.
If $\Phi_e(h) = 1\concat y\concat g$ for some $g\in\A - \C^-_{1,s}$ and $y\in\{0,1\}$
we also put $(1-y)\concat g$ into $\C_{1,s+1}$.

{\em Stage s+1=2e+2.} The construction to satisfy $R^1_e$ is completely symmetric
to the one for $R^0_e$, now using $\C_{1,s}$ instead of $\C_{0,s}$.
This ends the construction.

We verify that the construction succeeds in meeting all requirements.
At stage $s+1=2e+1$, the element $h$ put into $\C_0$ is a witness for
$\Phi_e(\C_0)\not\subseteq \B\meet \C_{1,s+1}$.
In order for $h$ to be a witness for $\Phi_e(\C_0)\not\subseteq\B\meet\C_1$
it suffices to prove that all elements $x\concat f$ entering $\C_1$ at a later
stage $t>2e+1$ are different from $\Phi_e(h)^-$.

If $\Phi_e(h)$ is not of the
form $1\concat y\concat g$ for $g\in\A-\C^-_{1,s}$ and $y\in\{0,1\}$ then this
is automatic, since only elements of this form are put into $\C_1$ at later stages.

Suppose $\Phi_e(h)$ is of the form $1\concat y\concat g$ for some
$g\in\A-\C^-_{1,s}$ and $y\in\{0,1\}$.
Then $(1-y)\concat g$ was put into $\C_{1,s+1}$ at stage $s+1$, if not earlier.
By construction, this ensures that all elements $x\concat f$ entering $\C_1$
at a later stage $t>s+1$ satisfy $f\neq g$:
\begin{itemize}

\item If $x\concat f$ enters $\C_{1,t+1}$ at $t= 2i+1$ then
$x\concat f = (1-y')\concat g'$ for some $g' \in \A-\C^-_{1,t}$ and $y'\in\{0,1\}$.
In particular $f\neq g$ since $g\in \C^-_{1,t}$.

\item If $x\concat f$ enters $\C_{1,t+1}$ at $t= 2i+2$ then
$f\in\A-\C^-_{1,t}$, so again $f\neq g$.

\end{itemize}
Thus $R^0_e$ is satisfied.
The verification of $R^1_e$ at stage $2e+2$ is again symmetric.
\end{proof}

\begin{lemma}  \label{genconstruction}
Let $\A$ and $\B$ be mass problems satisfying the condition 
\begin{center}
\hspace{\stretch{1}}
$\forall \C \subseteq \A \mbox{ finite } \; 
( \B\meet \C \not\leq_M \A )$. \hspace{\stretch{1}} 
\makebox[0pt][r]{\rm (\ref{condition})}
\end{center}
Then there exists an antichain $\C_\alpha$, $\alpha< 2^{\aleph_0}$,  
of mass problems such that $\C_\alpha \geq_M \A$ 
for every $\alpha$ and such that the elements 
$\B \meet \C_\alpha$ are also pairwise M-incomparable. 
\end{lemma}
\begin{proof}
We construct an antichain of size $2^{\aleph_0}$ as 
in Sacks' construction of such an antichain in the Turing degrees 
\cite[p462]{Odifreddi} by constructing a tree of $\C_\alpha$, $\alpha\in 2^\omega$, 
but now with the basic strategies from Lemma~\ref{construction}. 
As in Lemma~\ref{construction}, we build finite sets 
$\C_\sigma\subseteq \A\meet\A$, 
$\sigma\in 2^{<\omega}$. Given two sets $\C_\sigma$ and $\C_\tau$, 
$|\sigma| = |\tau| =s$, at stage $s=e$, we want to ensure that 
\begin{reqlist}

\item[$\Phi_e(\C_\beta)\not\subseteq \B\meet\C_\alpha$ \mbox{ and }]  
  
\item[$\Phi_e(\C_\alpha)\not\subseteq \B\meet\C_\beta$]

\end{reqlist}
for all $\alpha\sqsupset \sigma$ and $\beta\sqsupset \tau$. 
The basic strategy for doing this is exactly the same as in 
Lemma~\ref{construction}, and the way in which the strategies are 
put together on a tree is the same as in Sacks' construction. 
As a result, we obtain for every path $\alpha\in 2^\omega$ a set 
$\C_\alpha = \bigcup_{\sigma\sqsubset \alpha} \C_\sigma$ such that 
for every $e$ and $\beta\neq\alpha$ there is $f\in\C_\alpha$ such that 
$\Phi_e(f)\notin \B\meet \C_\beta$. 
So the sets $\B\meet \C_\alpha$, $\alpha\in 2^\omega$, are pairwise 
M-incomparable. 
\end{proof}

\begin{lemma}  \label{newconstruction}
Let $\A$ and $\B$ be mass problems satisfying the condition 
\begin{center}
\hspace{\stretch{1}}
$\forall \C \subseteq \A \mbox{ finite } \; 
( \B\meet \C \not\leq_M \A )$. \hspace{\stretch{1}} 
\makebox[0pt][r]{\rm (\ref{condition})}
\end{center}
Then there exists an antichain $\C_\alpha$, $\alpha< 2^{2^{\aleph_0}}$,  
of mass problems such that $\C_\alpha \geq_M \A$ 
for every $\alpha$ and such that the elements 
$\B \meet \C_\alpha$ are also pairwise M-incomparable. 
In particular, none of the $\C_\alpha$ is above $\B$.  
\end{lemma}
\begin{proof}
Start with the antichain $\C_{\alpha}$, $\alpha\in 2^\omega$, from 
Lemma~\ref{genconstruction}.  
If we knew that for every $\alpha$ there would be an $f\in \C_\alpha$ 
such that $f$ does not compute an element from any $\C_\beta$ with 
$\beta\neq\alpha$ then we 
could simply argue as in the original argument of Platek~\cite{Platek} 
showing that $\M$ has a big antichain, by taking $2^{2^{\aleph_0}}$ 
suitable combinations. But since we may not have this property we 
see ourselves forced to do something extra. 
For every $I\subseteq 2^\omega$ define 
$$
\C_I = 
\Meet_{\alpha\in I} \C_\alpha = 
\bigset{\alpha\oplus f: \alpha\in I \wedge f\in\C_\alpha}.
$$
Note that we use the indices explicitly to create a sort of disjoint 
union of the possibly continuum many $\C_\alpha$. This generalization 
of the meet operator to larger cardinalities than $\omega$ is no longer a 
natural meet operator, e.g.\ since the indices can be nontrivial now we loose the 
property that $\C_\alpha \geq_M \Meet_{\beta\in I} \C_\beta$, even if $\alpha\in I$, 
but this will not concern us. 

We want to construct a perfect set of indices $\TT\subseteq 2^\omega$ such that 
\begin{equation} \label{happy}
(\forall \alpha, \beta\in \TT)(\forall f\in \C_\alpha)(\forall g \in \C_\beta) 
\; [ \alpha\neq\beta \rightarrow \alpha\oplus f \T \beta\oplus g ].  
\end{equation}
The reason that it is possible to construct such a set of indices is that every 
$\C_\alpha$ is countable, and if $f\in \C_\alpha$ then $f$ in its totality is put 
into $\C_\alpha$ at some finite stage. 
We construct $\TT$ as the set of paths in a (noncomputable) tree 
$T\subseteq 2^{<\omega}$. 

{\em Construction of $T$.} 
Let $\C_\sigma$, $\sigma\in 2^{<\omega}$, refer to the finite approximations 
of the $\C_\alpha$ from the proof of Lemma~\ref{genconstruction}. 
At stage $s$ of the construction we have defined $T(\sigma)\in 2^{<\omega}$ 
for all $\sigma \in 2^{<\omega}$ of length $<s$. 
At stage $s$, for every $\sigma\neq \tau$ of length $s=e$ we guarantee 
{\setlength\arraycolsep{0pt}
\begin{eqnarray*} 
(\forall \alpha\sqsupset T(\sigma))(\forall \beta\sqsupset T(\tau))
(\forall f\in \C_\sigma)(\forall g \in \C_\tau) \; \big[ 
\Phi_e(\alpha\oplus f) && \neq \beta\oplus g  \wedge  \\
&& \Phi_e(\beta\oplus g) \neq \alpha\oplus f \big].
\end{eqnarray*}
}%
This can be realized in a standard finite extension construction \`{a} la Sacks, 
because the sets $\C_\sigma$ and $\C_\tau$ are finite. 
Given $f$ and $g$, the basic strategy for constructing $\alpha$ and $\beta$ with 
$\alpha\oplus f \T \beta\oplus g$ is the same as in the Kleene-Post 
construction of two incomparable sets. 
This concludes the construction of $T$.

The construction of $T$ guarantees that its set of paths 
$\TT$ satisfies~(\ref{happy}): 
Given $\alpha\neq \beta$ in $\TT$ and $f\in \C_\alpha$, $g\in \C_\beta$, 
the construction guarantees that 
$\Phi_e(\alpha\oplus f) \neq \beta\oplus g$ and   
$\Phi_e(\beta\oplus g) \neq \alpha\oplus f$ for all 
$e$ larger than the point in $T$ where $\alpha$ and $\beta$ split 
and larger than the stage where $f$ has entered $\C_\alpha$ and 
$g$ has entered~$\C_\beta$. 
Since we have this for almost every $e$, 
by padding we have $\alpha\oplus f \T \beta\oplus g$.

To finish the proof of the lemma, consider any family $\mathcal{I}$ of 
cardinality $2^{2^{\aleph_0}}$ of pairwise incomparable subsets of $\TT$ 
(cf.\ Proposition~\ref{antichains} below).  
If $I$ and $J$ are incomparable subsets of $\TT$ then by (\ref{happy}) 
we have that $\C_I \, |_M \, \C_J$. 
But then, since for every $\alpha$ we have $\C_\alpha\not\geq_M \B$, we 
also have 
$\B \meet \C_I \, |_M \, \B\meet \C_J$.
Note that $\C_I \geq_M \A$ because $\C_\alpha\subseteq \A\meet\A$. 
So the sets $\B\meet \C_I$, $I\in \mathcal{I}$, form an antichain of 
cardinality~$2^{2^{\aleph_0}}$. 
\end{proof}

\noindent
We note that one can prove the following
variant of Lemma~\ref{construction},
with a weaker hypothesis and a weaker conclusion, and with a
similar proof. 
A mass problem has {\em finite degree\/} if its M-degree contains 
a finite mass problem.

\begin{proposition} \label{variant} 
Let $\A$ be a mass problem that is not of finite degree, and
let $\B$ be any mass problem such that $\B\not\leq_M \A$.
Then there exists a pair $\C_0$, $\C_1$ of M-incomparable mass
problems above $\A$ such that neither of them is above $\B$.
\end{proposition}

An M-degree is a {\em degree of solvability\/} if it contains a
singleton mass problem. 
A mass problem is called {\em nonsolvable\/} if its M-degree
is not a degree of solvability. 
For every degree of solvability $\bfS$ there is a unique minimal M-degree 
$>\bfS$ that is denoted by $\bfS'$ (cf.\ Medvedev \cite{Medvedev}). 
If $\bfS = \deg_M(\{f\})$ then $\bfS'$ is the degree of the mass problem  
\begin{equation}  \label{Sprime}
\{f\}' = \bigset{n\concat g: f <_T g \wedge \Phi_n(g) = f }.
\end{equation}
(Note however that $\bfS'$ has little to do with the Turing jump.)
Dyment~\cite{Dyment} proved that the degrees of solvability are precisely 
characterized by the existence of such an~$\bfS'$. 
In particular the Turing degrees
form a first-order definable substructure of~$\M$.
The empty intervals in $\M$ are characterized by the following theorem. 
In view of what follows, it will be instructive to look at a proof of it. 

\begin{theorem} {\rm (Dyment \cite{Dyment})} 
\label{Dyment}
For Medvedev degrees $\bfA$ and $\bfB$ with $\bfA <_M \bfB$ it holds that
$(\bfA,\bfB)=\emptyset$ if and only if there is a degree of solvability $\bfS$
such that $\bfA \equiv_M \bfB \meet \bfS$, $\bfB\not\leq_M \bfS$, 
and $\bfB\leq_M \bfS'$.
\end{theorem}
\begin{proof}
(If) Suppose that $\bfS = \deg_M(\{f\})$ is as in the theorem and  
suppose that $\A\in \bfA$, $\B\in \bfB$, and 
$\B\meet \{f\}\leq_M \C \leq_M \B$. 
If $\C$ does not contain any element of Turing degree $\deg_T(f)$ 
then it follows that $\C \geq_M \B \meet \{f\}'$, because the elements 
of $\C$ that get sent to the $\{f\}$-side are all strictly 
above $f$, hence included in $\{f\}'$. So in this case $\C\geq_M \B$ 
by $\B\leq_M \{f\}'$.

Otherwise $\C$ contains an element of Turing degree $\deg_T(f)$, 
and consequently $\C\leq_M \{f\}$. 
Hence $\C \leq_M \B\meet \{f\} \equiv_M  \A$.

(Only if) Suppose that $(\A,\B)=\emptyset$. 
If $\A$ and $\B$ satisfy condition (\ref{condition}) then Lemma~\ref{construction} 
produces the M-incomparable sets $\B\meet\C_0$ and $\B\meet\C_1$ in 
$(\A,\B)$, so the interval is not empty in this case. 
So $\A$ and $\B$ do not satisfy condition (\ref{condition}) and  
hence there is a finite set $\C\subseteq \A$ such that 
$\B\meet \C\leq_M \A$. Then there is an $f\in\C$ such that 
$\{f\}\not\geq_M \B$, for otherwise we would have $\A\geq_M \B$.
Because the interval is empty we must have  $\A\equiv_M \B \meet \{f\}$ 
since there is no other possibility for $\B\meet \{f\}$. 
We also have $\B\meet \{f\}' \not\leq_M \A$ because 
both $\{f\}\not\geq_M \B$ and $\{f\}\not\geq_M \{f\}'$. 
Hence $\B\meet \{f\}' \geq_M \B$, again by emptiness of the 
interval, and in particular $\{f\}' \geq_M \B$. 
So we can take $\bfS$ to be $\deg_M(\{f\})$. 
\end{proof}

\begin{proposition}  \label{exactlytwo}
There are nonempty intervals in $\M$ that contain exactly 
two intermediate elements.
\end{proposition}
\begin{proof} 
Let $f$ and $g$ be T-incomparable and define 
$\A = \{f,g\}$, $\B= \{f\}'\meet\{g\}'$. 
We then have the situation as depicted in Figure~\ref{fig}. 
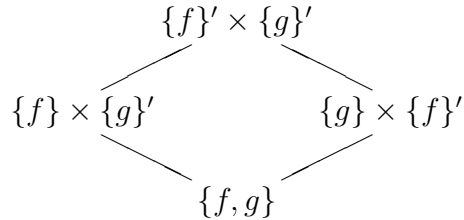
\begin{figure}[htb] 
\begin{center} 
\setlength{\unitlength}{1mm}
\begin{picture}(24, 29)
\put(12,0){\makebox[0cm][c]{$\{f, g\}$}} 
\put(-6,12){\makebox[-.5cm][c]{$\{f\}\meet \{g\}'$}}
\put(-6,16){\line(2,1){12}}   
\put(-6,10){\line(2,-1){12}}
\put(30,12){\makebox[.5cm][c]{$\{g\}\meet \{f\}'$}}
\put(30,16){\line(-2,1){12}}
\put(30,10){\line(-2,-1){12}}
\put(12,24){\makebox[0cm][c]{$\{f\}'\meet \{g\}'$}}
\end{picture}
\end{center}
\caption{An interval with exactly two intermediate elements.\label{fig}}
\end{figure} 
Note that $\{f\}\meet \{g\}'$ and $\{g\}\meet \{f\}'$ are indeed 
M-incomparable. 
By Theorem~\ref{Dyment} the two intervals 
$\big(\{f,g\}, \{f\}\meet \{g\}' \big)$ and 
$\big(\{f\}\meet \{g\}' , \{f\}'\meet\{g\}'\big)$ on the left side of 
the picture are empty, and by symmetry the same holds for the two 
intervals on the right side. 
Now suppose that 
\begin{equation} \label{RRR}
\{f,g\} \leq_M \C \leq_M \{f\}'\meet\{g\}'
\end{equation}
and that $\{f\}\meet \{g\}' \not\leq_M \C$. 
By Lemma~\ref{lemma:exactlyn} we then have $\C\leq_M \{g\}$. 
But then, since by (\ref{RRR}) we also have $\C\leq_M \{f\}'$, we have 
$\C\leq_M \{g\}\meet \{f\}'$. 
Thus if $\C \in (\A,\B)$ is not above 
$\{f\}\meet \{g\}'$ then it is below $\{g\}\meet \{f\}'$. 
Since all intervals depicted in Figure~\ref{fig} are empty, 
it follows that there are only the four possibilities 
for~$\C$. 
\end{proof}

\noindent
Since every interval in $\M$ is a lattice, 
we see that the interval of Figure~\ref{fig} is really isomorphic, 
as a lattice, to the Boolean algebra $\mathfrak{2}^2$. 
Proposition~\ref{exactlytwo} can be generalized to obtain finite intervals 
of size $2^n$, 
cf.\ Theorem~\ref{exactlyn}.  
We first prove a lemma. 

\begin{lemma} \label{lemma:exactlyn}
Let $n\geq 1$ and let $f_1,\ldots,f_n \in\omega^\omega$ be T-incomparable.
Suppose that $\C\geq_M \{f_1,\ldots,f_n\}$ and 
$\C\not\geq_M \{f_i\}' \meet \{f_j:j\neq i\}$. Then $\C\leq_M \{f_i\}$. 
\end{lemma}
\begin{proof} 
Let $\C$ satisfy the hypotheses of the lemma. 
By $\C\geq_M \{f_1,\ldots,f_n\}$ we have that  $\C$ is included in
the union of the Turing-upper cones of the $f_j$, $1\leq j\leq n$.
If $\C$ did not contain any element of degree $\deg_T(f_i)$ then
we would have $\{f_i\}' \meet \{f_j: j\neq i\}\leq_M \C$ as follows. 
Suppose that $\C\geq_M \{f_1,\ldots,f_n\}$ via $\Psi$. 
If $\Psi$ sends an element $h\in\C$ to some $f_j$, $j\neq i$, just let 
that happen, but if it sends $h$ to $f_i$ then instead output $h$.
We can recognize these distinctions effectively because
$\{f_1,\ldots,f_n\}$ is finite, so we can separate its elements 
by finite initial segments. 
This proves that $\C$ contains an element of degree $\deg_T(f_i)$, 
and consequently $\C\leq_M \{f_i\}$. 
\end{proof}

\begin{theorem} \label{exactlyn} 
Let $\B$ be any mass problem. 
Let $n\geq 1$ and let $f_1,\ldots,f_n \in\omega^\omega$ be T-incomparable 
such that $\{f_i\}\not\geq_M \B$ for every $i$. 
Then the interval 
$$
\big[\B\meet\{f_1,\ldots,f_n\}, \B\meet\{f_1\}'\meet\ldots\meet\{f_n\}'\big]
$$ 
is isomorphic to the Boolean algebra $\mathfrak{2}^n$. 
\end{theorem}
\begin{proof}
For $I\in \mathfrak{2}^n$ define 
$$
F(I) = \B\meet \Meet_{i\in I} \{f_i\}' \meet \{f_i: i\notin I\}.    
$$
Then clearly $F(I)\leq_M F(J)$ whenever $I\subseteq J$. 
Suppose that $I\neq J$, say $j\in J-I$. Then in $F(I)$ the factor $\{f_j\}$ 
occurs. But $\{f_j\}$ is neither above $\B$ nor above 
$\{f_j\}'$ nor above $\{f_i: i\neq j\}$, so $F(I)\not\geq_M F(J)$. 
So $F$ is an order-preserving injection. 
We verify that $F$ is onto. 
Suppose that $\C\in [F(\emptyset),F(n)]$.  
We prove that $\C$ is of the form $F(I)$ for some $I\subseteq n$. 
Let $I$ be a maximal subset of $n$ such that $\C\geq_M F(I)$. 
Note that such $I$ exists since $\C\geq_M F(\emptyset)$. 
Suppose that $i\notin I$ and that $\C$ contains no element of 
degree $\deg_T(f_i)$. 
Then similar argumentation as in Lemma~\ref{lemma:exactlyn} 
(just adding $\B$ to the argument) shows 
that $\C\geq_M F(I\cup\{i\})$, contradicting the maximality of $I$. 
So $\C$ contains an element of degree $\deg_T(f_i)$, 
and hence $\C\leq_M \{f_i\}$. 
Since we have this for every $i\notin I$ we have  
$\C\leq_M \{f_i: i\notin I\}$. 
Since we also have $\C\leq_M \B\meet \Meet_{i\in n} \{f_i\}'$ 
we have $\C\leq_M F(I)$.  Hence $\C \equiv_M F(I)$. 

We have proved that the interval is {\em order\/}-isomorphic 
to $\mathfrak{2}^n$. But then it follows automatically 
that it is isomorphic to $\mathfrak{2}^n$ as a {\em lattice}, 
since closing the elements $F(I)$ under $\meet$ and $\join$ 
cannot add any new elements because $F$ is onto. 
(It was already clear from the definition that the $F(I)$  
are closed under $\meet$.) 
Finally, since the interval is lattice-isomorphic to $\mathfrak{2}^n$, 
it follows that in fact it is a Boolean algebra itself. 
\end{proof}

\noindent
Platek~\cite{Platek} proved that $\M$ has the (for a collection  
of sets of reals maximal possible) cardinality $2^{2^{\aleph_0}}$ 
by showing that $\M$ has antichains of that cardinality. 
(The result was noted independently by Elisabeth Jockusch and John Stillwell.)
We now show that every interval in $\M$ is either of cardinality 
$2^n$ for some $n\geq 1$ or of cardinality $2^{2^{\aleph_0}}$.  
In particular, Theorem~\ref{exactlyn} is the {\em only\/} way 
to generate finite intervals. 
We will use the following lemma from~\cite{SorbiTerwijn}. 

\begin{lemma} {\rm (\cite{SorbiTerwijn})} \label{ok}
For any singleton mass problem $\S$,
if $\B\not\leq_M \S'$ then
$\S'$ and $\B$ satisfy condition~(\ref{condition}) from 
Lemma~\ref{newconstruction}.
\end{lemma}
\begin{proof}
Suppose that $\S = \{f\}$ and that $\C\subseteq \S'$ is finite such that
$\B \meet \C \leq_M \S'$, via $\Phi$ say. We prove that $\B\leq_M \S'$.

Recall the explicit definition of $\S'$ from equation~(\ref{Sprime}).
First we claim that for every $n\concat g \in \C$ there is $m\concat h \in \S'$
with $h\equiv_T g$ such that $\Phi(m\concat h)(0) = 0$,
that is, something from $\deg_T(g)$ is mapped to the $\B$-side.
To see this, let $m$ be such that $\Phi_m(f\oplus h') = f$ for all $h'$,
and let $h$ be of the form $f\oplus h'$ such that $\Phi(m\concat h)(0) =0$.
Such $h$ exists because $\C$ is finite, and for any number of finite elements
$\{f_0,\ldots,f_k\}$ strictly above $f$ it is always possible to build
$h >_T f$ such that $h$ is T-incomparable to all the $f_i$'s,
cf.\ \cite[p491]{Odifreddi}.
This $m\concat h$ is in $\S'$ and since it is incomparable to all 
the elements of $\C$ cannot be mapped by $\Phi$ to the $\C$-side, 
hence $\Phi(m\concat h)(0) = 0$. 
Now the computation $\Phi(m\concat h)(0) =0$ will use only a finite part of $h$,
so we can actually make $h$ of the same T-degree as $g$ by copying $g$ after
this finite part. This establishes the claim.

To finish the proof we note that from the claim it follows that $\B \leq_M \S'$:
If something is sent to the $\C$-side by $\Phi$ we can send it on to the
$\B$-side by the claim. Because $\C$ is finite we can do this uniformly.
More precisely, $\B \leq_M \S'$ by the following procedure.
By the claim fix for every $n\concat g \in \C$ a corresponding
$m\concat h \in \S'$ and a code $e$ such that $\Phi_e(g) = h$.
Given an input $n_0\concat g_0$, check whether $\Phi(n_0\concat g_0)(0)$
is $0$ or~$1$. (If it is undefined we do not have to do anything.) 
In the first case, output $\Phi(n_0\concat g_0)^-$,
i.e.\ $\Phi(n_0\concat g_0)$ minus the first element. This is then an element of $\B$.
In the second case $\Phi(n_0\concat g_0)^- \in \C$. Since $\C$ is finite we
can separate its elements by finite initial segments and determine
exactly which element of $\C$ $\Phi(n_0\concat g_0)^-$ is by inspecting only
a finite part of it. Now using the corresponding code $e$ that was chosen above
we output $\Phi\big(m\concat \Phi_e\big( \Phi(n_0\concat g_0)^-\big)\big)$,
which is again an element of $\B$.
\end{proof}

\begin{theorem} \label{dichotomy} 
Let $[\bfA,\bfB]$ be an interval in $\M$ with $\bfA <_M \bfB$. 
Then either 
$[\bfA,\bfB]$ is isomorphic to the Boolean algebra $\mathfrak{2}^n$ 
for some $n\geq 1$, 
or $[\bfA,\bfB]$ contains an antichain of size $2^{2^{\aleph_0}}$.  
\end{theorem}
\begin{proof}
Let $\A$ and $\B$ be mass problems of degree $\bfA$ and $\bfB$, respectively. 
If $\A$ and $\B$ satisfy condition (\ref{condition})    
then Lemma~\ref{newconstruction} immediately gives an antichain of  
size $2^{2^{\aleph_0}}$ between $\A$ and $\B$.   
Suppose next that $\A$ and $\B$ do not satisfy condition (\ref{condition}): 
Let $\C\subseteq \A$ be finite such that $\B\meet\C \leq_M \A$. Since 
also $\A\leq_M \B\meet\C$ we then have $\A \equiv_M \B\meet \C$. 
Since $\C$ is finite we can separate its elements by finite initial segments 
and hence it holds that 
$$
\B\meet \C \equiv_M \B\meet \big\{f\in \C: \{f\}\not\geq_M \B \wedge 
f \mbox{ is of minimal T-degree in } \C \big\}, 
$$ 
so we may assume without loss of generality that the 
elements of $\C$ are pairwise T-incomparable and satisfy $\{f\}\not\geq_M \B$.

If $\C=\emptyset$ then $\A \equiv_M \B$ so the interval contains 
just this one element. 

Suppose that $\C = \{f_1,\ldots, f_n\}$, $n\geq 1$, so that 
$\A \equiv_M \B \meet \{f_1,\ldots, f_n\}$. 
If $\B \leq_M \{f_1\}'\meet\ldots\meet \{f_n\}'$ then 
$\B \equiv_M \B\meet \{f_1\}'\meet\ldots\meet \{f_n\}'$, so by 
Theorem~\ref{exactlyn} 
$[\A,\B]$ is isomorphic to $\mathfrak{2}^n$. 

If $\B \not\leq_M \{f_1\}'\meet\ldots\meet \{f_n\}'$ then 
there is an $i$ such that $\B\not\leq_M \{f_i\}'$. 
We now apply Lemma~\ref{newconstruction} to $\{f_i\}'$ and $\B$.
This is possible because $\{f_i\}'$ and $\B$ satisfy
condition~(\ref{condition}) by Lemma~\ref{ok}.
Lemma~\ref{newconstruction} now produces an antichain of elements 
$\B\meet\C_\alpha$ with $\C_\alpha\geq_M \{f_i\}'$. 
The elements of the antichain are clearly below $\B$, 
and they are also above $\A$ since
$\C_\alpha \geq_M \{f_i\}' \geq_M \{f_i\}\geq_M \A$.
So we have again an antichain of
size $2^{2^{\aleph_0}}$in the interval $(\A,\B)$.
\end{proof}

\begin{corollary} {\rm (Sorbi and Terwijn \cite{SorbiTerwijn})}\label{nonlinear}
If $(\bfA,\bfB)\neq\emptyset$ then there is a pair of incomparable 
degrees in $(\bfA,\bfB)$. 
\end{corollary}


\noindent
In \cite{SorbiTerwijn} Corollary~\ref{nonlinear} was used to show that 
the linearity axiom
$$
(p\rightarrow q)\vee(q\rightarrow p)
$$ 
is not in any of the theories $\Th(\M/\A)$ for $\A>_M 0'$, 
where $\Th(\M/\A)$ is the set of all propositional formulas that are 
valid on the Brouwer algebra $\M/\A$. 

Note that there are both $2^{2^{\aleph_0}}$ examples of finite and of  
infinite intervals in~$\M$. 
For the first, note that if $\B$ is upwards closed under $\leq_T$ and 
$f \notin \B$ then by Theorem~\ref{exactlyn} we can associate 
a finite interval with the pair $(\B,f)$. 
Now as in Plateks argument, taking an antichain of size $2^{\aleph_0}$ 
in the Turing degrees we see that there are $2^{2^{\aleph_0}}$ 
such pairs, all defining different finite intervals. 

To see that there are also $2^{2^{\aleph_0}}$ infinite intervals, 
note that by the proof of Theorem~\ref{dichotomy} it suffices to 
show that there are $2^{2^{\aleph_0}}$ pairs $(\A,\B)$ 
(with all the $\A$'s of different M-degree) satisfying 
condition~(\ref{condition}). But this is easy to see, 
again using the antichain from above.

\section{Notes about chains and antichains in $\P(\kappa)$} \label{notes}

As a preparation for the next section we collect  some notes 
about chains and antichains in $\mbox{}^\kappa 2$, for an arbitrary 
cardinal $\kappa$. 
We claim no originality, but include two simple proofs 
for later reference. 
Our set-theoretic notation follows Kunen~\cite{Kunen}. 
Note however that by an antichain we just mean a set of pairwise 
incomparable elements, whereas Kunen uses a stronger notion 
(with ``incompatible'' instead of ``incomparable''). 
A chain in $\P(\kappa)$ is any family of subsets of $\kappa$ that is 
strictly linearly ordered by $\subsetneq$.

\begin{proposition} \label{antichains}
Let $\kappa$ be any cardinal. Then the partial order $(\P(\kappa),\subseteq)$ 
has an antichain of size $2^\kappa$.
\end{proposition}
\begin{proof}
By explicit construction. 
Build a tree of $C_X\subseteq \kappa$, $X\in \mbox{}^\kappa 2$, 
in $\kappa$ stages as follows. 
Start with $C_0 = \emptyset$. 
For any stage $\alpha<\kappa$, pick two fresh (i.e.\ not previously used 
in the construction) elements of $\kappa$, and 
for any $\sigma \in \mbox{}^{<\kappa}2$ of length $\alpha$, 
put one fresh element in $C_{\sigma 0}$ and the other in~$C_{\sigma 1}$. 
At limit stages $\alpha$ take unions, that is, let 
$C_\sigma = \bigcup_{\tau\sqsubset\sigma} C_\tau$ for  
any $\sigma\in \mbox{}^{<\kappa}2$ of length $\alpha$. 
For every $X\in \mbox{}^\kappa 2$ 
let $C_X = \bigcup_{\alpha< \kappa} C_{X\smallrestr \alpha}$.  
End of construction. 

Clearly, if $X\neq X'$ then $C_X | C_{X'}$ 
(look at the first $\alpha$ where $X$ and $X'$ differ), 
so the $C_X\subseteq \kappa$ form an antichain of size $2^\kappa$. 
\end{proof}

\noindent
So for antichains we do not need any special properties of $\kappa$. 
For chains the situation is more complicated, but we have the following 
result. 

\begin{proposition} \label{chains}
Let $\kappa$ be any cardinal with $2^{<\kappa}=\kappa$. 
Then the partial order $(\P(\kappa),\subseteq)$ has a chain of size $2^\kappa$. 
\end{proposition}
\begin{proof}
This is a generalization of the fact that $\P(\omega)$ has chains of 
size $2^\omega$. 
We view $\P(\kappa)$ as the set of paths in the tree $\mbox{}^{<\kappa}2$. 
Let $<_L$ denote the Kleene-Brouwer ordering on 
$\mbox{}^{<\kappa}2 \cup \mbox{}^\kappa 2$:   
for strings $\sigma$ and $\tau$ of length $\leq \kappa$, 
$\sigma <_L \tau$ if there is an $\alpha<\kappa$ such that 
$\sigma(\alpha) < \tau(\alpha)$, 
i.e.\ $\sigma$ branches off to the left of $\tau$ in the tree $\mbox{}^{<\kappa}2$.   
Now every $X\in \mbox{}^\kappa 2$ has an associated ``Dedekind cut'' 
$$
{\rm Cut}(X) = \bigset{\sigma\in \mbox{}^{<\kappa}2: \sigma<_L X }. 
$$
By $2^{<\kappa}=\kappa$, every cut corresponds to a subset of $\kappa$, 
and we have that $X<_L Y$ implies that ${\rm Cut}(X)\subsetneq {\rm Cut}(Y)$. 
So under this assumption the cuts form a chain of size $2^\kappa$.
\end{proof}

\noindent
Note that the chains in Proposition~\ref{chains} cannot be {\em well-ordered\/}, 
since well-ordered chains in $\P(\kappa)$ have size at most $\kappa$ 
(since every next element of the chain has to add a new element). 

In Section~\ref{ChainsinM} we will be interested in the case $\kappa = 2^\omega$, 
the Medvedev lattice being a collection of sets of reals. 
The fact that for $\kappa = 2^\omega$ the condition 
$2^{<\kappa}=\kappa$ of Proposition~\ref{chains} is independent 
of ZFC (since it is true under CH and false e.g.\ when  
$2^\omega = \omega_2$ and $2^{\omega_1} = \omega_3$) 
suggests that the existence of big chains in $\P(2^\omega)$  
might also be independent. 

\begin{problem} \label{problem}
Settle the independence of ZFC of the existence of chains of size 
$2^{2^{\aleph_0}}$ in $\P(2^\omega)$. 
\end{problem}

\noindent
Problem~\ref{problem} was put before several set theorists. 
Perhaps surprisingly, it seems that it is open.

\section{Chains in $\M$} \label{ChainsinM}

Since $\M$ is defined by factoring out sets of reals modulo the 
reduction relation $\leq_M$, a priori its maximal possible cardinality 
is $2^{2^{\aleph_0}}$. 
In Section~\ref{intervals} we have seen that $\M$ has indeed cardinality   
$2^{2^{\aleph_0}}$, and in fact that every infinite interval contains an 
antichain of this cardinality.  
In this section we consider the height of intervals, that is, 
we discuss chains. First we note that every M-degree is as large as 
set-theoretically possible: 

\begin{proposition} \label{degrees}
For every mass problem $\F\subseteq 2^\omega$ we have 
$|\deg_M(\F)|= 2^{2^{\aleph_0}}$. 
\end{proposition}
\begin{proof}
This follows by simple counting. For every $\X\subseteq 2^\omega$ define 
$$
\A_\X = \bigset{f\oplus g : f\in \F \wedge g \in\X\cup\{0^\omega\} }, 
$$
where $0^\omega$ is the all zero sequence. 
Then $\A_\X \leq_M \F$ since for all $f\in\F$, $f\oplus 0^\omega\in \A_\X$, 
and 
$\F\leq_M \A_\X$ since for all $h\in\A_\X$, $h_0 =f \in\F$, where 
$h_0$ is the unique component such that $h = h_0 \oplus h_1$. 
So $\A_\X \equiv_M \F$ for every $\X$, hence the result follows. 
\end{proof}

\noindent
We have seen in Section~\ref{notes} that whether or not $\P(2^\omega)$ has 
chains of size $2^{2^{\aleph_0}}$ may depend on set-theoretic properties of 
$2^\omega$. 
The same holds for $\M$. 
Next we show that it is at least consistent with ZFC that $\M$ has chains 
of the size of its own cardinality. 

\begin{theorem}  \label{chainsM}
CH implies that $\M$ has a chain of size $2^{2^{\aleph_0}}$. 
\end{theorem}
\begin{proof}
We build on the proof of Proposition~\ref{chains}. 
Note that under CH we have $2^{<2^\omega} = 2^{<\omega_1} = 2^\omega$, 
so the condition of Proposition~\ref{chains} is satisfied for $\kappa = 2^\omega$.   
Let $\{f_\alpha: \alpha< 2^\omega\}$ be a set of pairwise Turing incomparable   
elements of $\mbox{}^\omega 2$, which exists by a result of Sacks,   
cf.\ \cite{Odifreddi}. 
Now by CH the sets ${\rm Cut}(X)$ defined in the proof of Proposition~\ref{chains} 
correspond with subsets of $\mbox{}^\omega 2$, hence also with subsets of 
$\{f_\alpha: \alpha< 2^\omega\}$. Call this correspondence $F$, so that 
${\rm Cut}(X)$ corresponds to $F({\rm Cut}(X))$. 
Again we have $X<_L Y$ implies that ${\rm Cut}(X) \subsetneq {\rm Cut}(Y)$, 
which in turn implies that $F({\rm Cut}(X)) \subsetneq F({\rm Cut}(Y))$.  
But then $F({\rm Cut}(Y)) <_M F({\rm Cut}(X))$ because any $f_\alpha$ 
in $F({\rm Cut}(Y))- F({\rm Cut}(X))$ cannot compute an element of 
$F({\rm Cut}(X))$. 
\end{proof}

\noindent
From the proof of Theorem~\ref{chainsM} we see that the conditions for the 
existence of big chains in $\P(2^\omega)$ and in $\M$ are exactly  the same.

In Theorem~\ref{dichotomy} we saw that every interval in $\M$ is either 
small (isomorphic to a finite Boolean algebra) or 
contains a big antichain. We now show that it is consistent that 
the interval also has a big chain whenever it has a big antichain. 

\begin{theorem} \label{extendeddichotomy} 
Let $[\bfA,\bfB]$ be an interval in $\M$ with $\bfA <_M \bfB$. 
Then either 
$[\bfA,\bfB]$ is isomorphic to the Boolean algebra $\mathfrak{2}^n$ 
for some $n\geq 1$, 
or $[\bfA,\bfB]$ contains an antichain of size $2^{2^{\aleph_0}}$. 
In the latter case, assuming CH, it also contains 
a chain of size $2^{2^{\aleph_0}}$. 
\end{theorem}
\begin{proof} By Theorem~\ref{dichotomy} and its proof we only have to 
show how to obtain a big chain from the big antichain we constructed in 
Lemma~\ref{newconstruction}. 
Recall the sets $\C_\alpha$ and $\C_I$ from the proof of 
Lemma~\ref{newconstruction}, and also the special set of indices 
$\TT\subseteq 2^\omega$ defined there. 
$\TT$ is of cardinality $2^{\aleph_0}$ and order-isomorphic to $2^\omega$, 
with the order inherited from $2^\omega$. 
Let $\leq_L$ be the Kleene-Brouwer ordering defined 
as in Proposition~\ref{chains}, but now on the 
tree $\mbox{}^{< \TT}2 \cup \mbox{}^{\TT}2$. 
For every $I\in \mbox{}^\TT 2$ we have the associated set 
$$
{\rm Cut}(I) = \bigset{\sigma\in \mbox{}^{<\TT}2: \sigma\leq_L I }. 
$$
By CH, $|\mbox{}^{<\TT}2| = 2^{<2^\omega} = 2^\omega = |\TT|$, so we can 
associate with every ${\rm Cut}(I)$ a subset $F({\rm Cut}(I))$ of $\TT$. 
Now let 
$$
\mathcal{E}(I) = \bigset{\alpha\oplus f: 
\alpha \in F({\rm Cut}(I)) \wedge f\in \C_\alpha}.   
$$ 
Clearly the cuts, and hence the $\mathcal{E}(I)$, form a monotone sequence, 
that is, 
$$
J \leq_L I \Rightarrow {\rm Cut}(J) \subseteq {\rm Cut}(I) 
\Rightarrow \mathcal{E}(J) \subseteq \mathcal{E}(I) 
\Rightarrow \mathcal{E}(I) \leq_M \mathcal{E}(J).
$$
The sequence is strict because $J <_L I$ implies that there is an 
$\alpha\in F({\rm Cut}(I)) - F({\rm Cut}(J))$, 
and by the property~(\ref{happy}) in the proof of Lemma~\ref{newconstruction} we then have 
that $\mathcal{E}(I) \not\geq_M \mathcal{E}(J)$. 
Thus we have  $\mathcal{E}(I) <_M \mathcal{E}(J)$ whenever $J<_L I$. 
So the sets $\mathcal{E}(I)$ form a chain in $\M$  of cardinality 
$|\mbox{}^\TT 2| = 2^{2^{\aleph_0}}$. 
\end{proof}

\section{Embeddings into $\M$}   \label{embeddings}

Sorbi characterized the countable lattices that are embeddable into $\M$ 
as follows: 

\begin{theorem} {\rm (Sorbi~\cite{Sorbi1990, Sorbi})} \label{embedding}
A countable distributive lattice with 0,1 is embeddable
into $\M$ (preserving 0 and 1) 
if and only if 0 is meet-irreducible and
1 is join-irreducible.
\end{theorem}

\noindent
Sorbi proved Theorem~\ref{embedding} by embedding the (unique) 
countable dense Boolean algebra into $\M$. 
Below we show that this embedding is optimal as far as 
cardinalities are concerned: Every Boolean algebra embeddable into $\M$ 
must be countable. We first show that the (dual of) the 
large Boolean algebra $\P(2^\omega)$ is embeddable into $\M$ as 
an upper semilattice, i.e.\ preserving joins but not necessarily meets. 

Let $f =^* g$ denote that the functions $f$ and $g$ differ only on finitely 
many elements. 

\begin{lemma} \label{maart2006}
There exist a noncomputable $g\in 2^\omega$ and $g_X\in 2^\omega$, $X\in 2^\omega$, 
and a 
computable functional $\Psi$ such that for all $X$ and~$Y$, 
\begin{itemize}
\item $g \not\leq_T g_X$, 
\item $X\neq Y \Longrightarrow g_X \T g_Y$, 
\item $X\neq Y \Longrightarrow \Psi(g_X\oplus g_Y) =^* g$, 
\item $g_X \leq_T g \oplus X$. 
\end{itemize}
\end{lemma}
\begin{proof}
Note that this is an extension of the existence of 
an antichain of size $2^{\aleph_0}$ in the Turing degrees \cite[p462]{Odifreddi}. 
It can be realized by standard methods from computability theory, so 
we only sketch the idea. 
Construct a c.e.\ set $G$ and a tree of finite sets $G_\sigma$, 
$\sigma\in 2^{<\omega}$, such that for the sets 
$G_X = \bigcup_{\sigma\sqsubset X} G_\sigma$, with $X\in 2^\omega$, we have that 
$X\neq Y \Longrightarrow G_X \T G_Y$ and $G_X \cup G_Y =^* G$. 
Of course we cannot make all the $G_X$ c.e.\ because there are too many 
of them, but at least we can make them all c.e.\ relative to the path $X$ 
that defines them. 
Given a string $\sigma\in 2^{<\omega}$, to make all the paths extending 
$\sigma0$ T-incomparable with those extending $\sigma1$ we have standard 
Friedberg-Muchnik requirements 
$$
G_{\sigma0} \neq \{e\}^{G_{\sigma1}}, \;\; G_{\sigma1} \neq \{e\}^{G_{\sigma0}},
$$
that can be satisfied using the usual strategy and that can be put together 
in a finite injury argument. Every element enumerated into some $G_\sigma$ 
is simultaneously enumerated into $G$. 
If an element $x$ is enumerated at stage $s$, we enumerate it in 
{\em all\/} $G_{\sigma}$ with $|\sigma|=s$, except one. 
This may cause injuries but these will be finitary. 
It ensures that if $X\neq Y$, then from some stage onwards all 
$x$ entering $G$ enter at least one of $G_X$ and $G_Y$, hence 
$G_X \cup G_Y =^* G$. Since the construction is effective, 
to decide which elements are in $G_X$ it suffices 
to know $G$ and the path $X$, hence $G_X\leq_T G\oplus X$.  
It now follows that $G \not\leq_T G_X$ for all $X$: 
Suppose that $G\leq_T G_X$. Let $Y$ be computable and different from~$X$. 
We then have $G_Y \leq_T G \oplus Y \leq_T G_X$, contradicting 
that $G_X \T G_Y$. 
From this the lemma clearly follows. 
\end{proof}

\noindent
For any Boolean algebra $\B$, let $\dual{\B}$ denote the dual of $\B$. 

\begin{theorem} \label{largeupper}
There is an embedding of $\dual{\P(2^\omega)}$ into $\M$ as an upper semilattice. 
\end{theorem}
\begin{proof}
Let $g$ and $g_X$ be as in Lemma~\ref{maart2006}. 
Let ${\rm Fin(g)}$ be the set of all finite differences of $g$. 
Consider the mass problems 
$$
\A_I = {\rm Fin}(g) \meet \{g_X: X\in I\}
$$
for every $I\subseteq 2^\omega$. 
We claim that $F:\dual{\P(2^\omega)} \hookrightarrow \M$ defined by 
$I\mapsto \deg_M(\A_I)$ is an embedding of upper semilattices. 
Clearly $I\subseteq J$ implies that $F(J) \leq_M F(I)$, and we have 
$F(J) <_M F(I)$ if the inclusion is strict because $g_X$ can neither 
compute $g$ nor any other $g_Y$. 
We check that $F(I\cap J) \equiv_M F(I) \join  F(J)$. 
Clearly $\A_I$, $\A_J\leq_M \A_{I\cap J}$ via inclusion, hence 
$\A_I \join \A_J \leq_M \A_{I\cap J}$. 
Conversely, $\A_{I\cap J} \leq_M \A_I \join \A_J$:  
Suppose $g_i \in \A_I$ and $g_j \in \A_J$. 
If one of the two is in $0\,\concat\,{\rm Fin}(g)$ (which we can see from the first 
bit) then we are immediately done. 
Otherwise, and if $g_i\neq g_j$, it holds that 
$\Psi(g_i^-\oplus g_j^-) =^* g$ by Lemma~\ref{maart2006}, so again 
we can produce an element of $0\,\concat\,{\rm Fin}(g)$, hence of $\A_{I\cap J}$. 
If $g_i = g_j$ then $g_i \in \A_{I\cap J}$, so in this case we could just output 
$g_i$. However, we cannot a priori distinguish between this case and the previous 
one, so given $g_i$ and $g_j$ both not in $0\,\concat\,{\rm Fin}(g)$ we start 
outputting $g_i$ until, if ever, a difference with $g_j$ is found, in which 
case we continue outputting $\Psi(g_i^-\oplus g_j^-)$. 
In the latter case we will output a finite difference of a finite 
difference of $g$, which is in $\A_{I \cap J}$. 
\end{proof}

\noindent
The next result shows that the Boolean algebra 
$\dual{\P(2^\omega)}$ is not embeddable into $\M$ as a lattice. 
Note that such an embedding would automatically be 
an embedding as a Boolean algebra.  

\begin{theorem} \label{countable}
Suppose $\B$ is a Boolean algebra that is embeddable into $\M$ as a lattice 
(i.e.\ preserving meets and joins). Then $\B$ is countable. 
\end{theorem}
\begin{proof}
Suppose for a contradiction that $\B$ is uncountable and that 
$F:\B\hookrightarrow \P(\omega^\omega)$ defines an embedding of $\B$, 
i.e.\ that $X\mapsto \deg_M(F(X))$ is an embedding of $\B$ into~$\M$. 
By the Stone representation theorem we may think of $\B$ as an 
algebra of sets, so we denote the lattice operations in $\B$ by 
$\cap$ and $\cup$, 
let $\emptyset$ be the bottom element of $\B$, and for $X\in \B$ let 
$\cmp{X}$ denote the complement of $X$ in $\B$.  
By assumption we then have 
$F(X\cap Y) \equiv_M F(X) \meet F(Y)$ and 
$F(X\cup Y) \equiv_M F(X) \join F(Y)$ for all $X$ and $Y$ in~$\B$. 
In particular we have for every $X\in\B$ that 
\begin{equation} \label{many}
F(X) \meet F(\cmp{X}) \leq_M F(\emptyset). 
\end{equation}
Because $\B$ is uncountable there are uncountably many inequalities of the 
form~(\ref{many}). 
Since there are only countably many computable functionals $\Phi_e$, 
there must be two different $X$, $Y\in B$ such that 
$F(X) \meet F(\cmp{X}) \leq_M F(\emptyset)$ and 
$F(Y) \meet F(\cmp{Y}) \leq_M F(\emptyset)$ 
via the {\em same\/}~$\Phi$. 
Since $X\neq Y$ we have that either $X\not\subseteq Y$ or 
$Y\not\subseteq X$. For definiteness say that $X\not\subseteq Y$. 
We have  
$$
F(\cmp{X})\cap F(\cmp{Y}) \geq_M F(\cmp{X}) \join F(\cmp{Y}) 
\equiv_M F(\cmp{X} \cup \cmp{Y}), 
$$
and therefore    
$$  
\begin{array}{rcll}
F(\emptyset) & \geq_M & F(X) \meet 
\big(F(\cmp{X})\cap F(\cmp{Y})\big) & \mbox{(via $\Phi$)} \\
& \geq_M & F(X) \meet F(\cmp{X} \cup \cmp{Y}) &  \\
&\equiv_M& F(X \cap (\cmp{X}\cup \cmp{Y})) & \\
&\equiv_M& F(X\cap \cmp{Y}) & \\
&\geq_M& F(\emptyset). & 
\end{array}
$$  
Hence $F(\emptyset) \equiv_M F(X\cap \cmp{Y})$, which is a contradiction 
since $X\cap \cmp{Y} \neq \emptyset$ because $X\not\subseteq Y$. 
\end{proof}

\section{The Muchnik lattice}  

There is a nonuniform variant of the Medvedev lattice, called 
the Muchnik lattice, that was introduced by Muchnik in \cite{Muchnik}. 
This is the structure $\M_w$ resulting from the reduction relation on 
mass problems defined by 
$$
\A \leq_w \B \equiv  (\forall f\in \B)(\exists g\in \A)[g\leq_T f]. 
$$
That is, every solution to the mass problem $\B$ can compute a solution 
to the mass problem $\A$, but maybe not in a uniform way. 
$\M_w$ is a distributive lattice in the same way that $\M$ is, with 
the same lattice operations and $\bf 0$ and~$\bf 1$. 
An M-degree is a {\em Muchnik degree\/} if it contains 
a mass problem that is upwards closed under Turing reducibility 
$\leq_T$. The Muchnik degrees of $\M$ form a substructure that is 
isomorphic to $\M_w$. 

We check which of the results from the previous sections hold also for 
$\M_w$ instead of $\M$, replacing $\leq_M$ by $\leq_w$.  
As we will see, because of the lack of uniformity much more is 
possible in $\M_w$, making it a structure much closer to 
the Turing degrees. In particular we no longer have strong dichotomies as 
in Theorem~\ref{dichotomy}. 

\begin{example}  \label{3chain}  
Let $f$ and $\emptyset <_T g <_T f$ be such that the Turing lower cone of $f$ 
consists precisely of three elements: 
$$
\forall h \big(h\leq_T f \rightarrow h \mbox{ computable } \vee 
h \equiv_T g \vee h \equiv_T f \big).
$$ 
Such $f$ exists since Titgemeyer proved that the three-element chain is 
embeddable into the Turing degrees as an initial segment, 
cf.\ \cite[p526]{Odifreddi}.   
Now let $\B 
= \bigset{h : h \not\leq_T f}$ and $\A = \B \meet \{g\}$. 
We claim that $(\A,\B)$ contains only 
one element, namely $\B\meet \{f\}$. Suppose that $\C\in(\A,\B)$. 
Then there exists $h \in \C$ such that $\{h\}\not\geq_w \B$, 
hence $h\leq_T f$, so either $h\equiv_T g$ or $h\equiv_T f$. 
If $\C$ contains such an $h$ with $h\equiv_T g$ then 
$\C \leq_w \B \meet \{g\} = \A$. 
Otherwise, all $h\in C$ with $\{h\}\not\geq_w \B$ have $h\equiv_T f$, 
so we have both $\C\leq_w \B\meet \{f\}$ and $\B \meet \{f\} \leq_w \C$. 
\end{example}

\noindent
From Example~\ref{3chain} we see that there are linear nonempty intervals in 
$\M_w$. This shows in particular that Corollary~\ref{nonlinear} fails for $\M_w$. 

In the proof of Dyments Theorem~\ref{Dyment} given in Section~\ref{intervals} 
we used Lemma~\ref{construction}, which fails for $\M_w$  
(because $\M_w$ contains an infinite downward chain, 
cf.\ the proof of Lemma~\ref{subst}). 
Nevertheless, the theorem still holds for $\M_w$. 
Instead of Lemma~\ref{construction} one can use the following much easier result:  

\begin{lemma} \label{subst}
Suppose that $\A$ and $\B$ satisfy 
\begin{equation}   \label{new}
\forall \C \subseteq \A \mbox{ finite } \; ( \B\meet \C \not\leq_w \A ). 
\end{equation}
Then there exists $\C \geq_w \A$ such that $\C \not\geq_w \B$ and 
$\B\meet \C \not\leq_w \A$.  
If moreover $\A \leq_w \B$ then the interval $(\A, \B)$ is infinite. 
\end{lemma}
\begin{proof} 
Since from (\ref{new}) it follows that $\B\not\leq_w \A$, there is 
$f\in \A$ such that $\{f\} \not\geq_w \B$. 
Again by (\ref{new}) we have that $\B\meet \{f\} \not\leq_w \A$, 
so we can take $\C = \{f\}$. 
 
If in addition $\A \leq_w \B$ then we have 
$\A <_w \B \meet \{f\} <_w \B$. 
Since $\A$ and $\B \meet \{f\}$ also satisfy (\ref{new}) we can 
by iteration of the first part of the lemma obtain an infinite 
downward chain in $(\A, \B)$. 
\end{proof}

\noindent 
The proof of Theorem~\ref{Dyment} for $\M_w$ is now exactly the same 
as the proof given above, replacing $\leq_M$ by $\leq_w$ and 
using Lemma~\ref{subst} where previously Lemma~\ref{construction} was used. 

Theorem~\ref{exactlyn} still holds for $\M_w$, with the same proof, 
but as we have seen in Example~\ref{3chain} it is not the only way anymore 
to generate finite intervals. %
%
%
%
In Terwijn~\cite{Terwijnta} the finite intervals of $\M_w$ are characterized 
as a certain proper subclass of the finite distributive lattices. 

In contrast to Theorem~\ref{extendeddichotomy}, 
an interval in $\M_w$ can be infinite without containing a large antichain.
In fact, in \cite{Terwijnta} it is proved that 
there are intervals in $\M_w$ with maximal antichains of every possible size. 
Similar results can be obtained for chains. 
Theorem~\ref{chainsM} holds for $\M_w$, with the same proof,
so the conditions for the existence of chains of size $2^{2^{\aleph_0}}$
in $\M$ and in $\M_w$ are the same.
The consistency of the existence of chains of this size also
follows from Proposition~\ref{embeddingMw} below.

Proposition~\ref{degrees} also holds for $\M_w$, with the same proof, 
so again every Muchnik degree is as large as set-theoretically possible. 
Theorem~\ref{countable} does not hold for $\M_w$, as we can in fact 
embed the dual of $\P(2^\omega)$: 

\begin{proposition} \label{embeddingMw}
$\dual{\P(2^\omega)}$ is embeddable into $\M_w$ as a Boolean algebra. 
\end{proposition}
\begin{proof}
This follows simply by noting that the embedding $F$ given in the proof 
of Theorem~\ref{largeupper}, that did not preserve meets for $\M$, 
does in fact preserve meets in $\M_w$. 
For $\M$ we have  
$F(I\cup J) \leq_M F(I) \meet F(J)$ but not necessarily 
$F(I\cup J) \geq_M F(I) \meet F(J)$. 
But we do have $F(I\cup J) \geq_w F(I) \meet F(J)$, as is immediate 
from the definitions. 
\end{proof}

\begin{acknowledgments}
The author thanks Martin Goldstern, Jakob Kellner, and Andrea Sorbi 
for helpful discussions. 
\end{acknowledgments}

\end{document}